\documentclass[12pt]{amsart}
\usepackage{latexsym}

\newtheorem{thm}{Theorem}[section]
\newtheorem{propn}[thm]{Proposition}
\newtheorem{lemma}[thm]{Lemma}
\newtheorem{coro}[thm]{Corollary}
\newtheorem{rem}[thm]{Remark}

\begin{document}


\title{Deformations of the generalised Picard bundle}

\author{I. Biswas}

\address{School of Mathematics, Tata Institute of
Fundamental Research, Homi Bhabha Road, Bombay 400005, India}

\email{indranil@math.tifr.res.in}

\author{L. Brambila-Paz}

\address{CIMAT, Apdo. Postal 402, C.P. 36240. Guanajuato, Gto,
M\'exico}

\email{lebp@cimat.mx}

\author{P. E. Newstead}

\address{Department of Mathematical Sciences, The University of
Liverpool, Peach Street, Liverpool, L69 7ZL, England}

\email{newstead@liverpool.ac.uk}


\thanks{All authors are members of the research group VBAC
(Vector Bundles on Algebraic Curves), which is partially supported
by EAGER (EC FP5 Contract no. HPRN-CT-2000-00099) and by EDGE (EC
FP5 Contract no. HPRN-CT-2000-00101).  The second author acknowledges
the support of CONACYT grant 40815-F}

\date{28/08/03}

\begin{abstract}

Let $X$ be a nonsingular algebraic curve of genus $g\geq 3$, and let
${\mathcal M}_{\xi}$ denote the moduli space of stable
vector bundles of rank $n\ge2$ and degree $d$ with fixed determinant
$\xi$ over $X$ such that $n$ and $d$ are
coprime. We assume that if $g=3$ then $n\geq 4$
and if $g=4$ then $n\geq 3$, and suppose further that $n_0$, $d_0$ are integers
such that $n_0\ge1$ and $nd_0+n_0d>nn_0(2g-2)$. Let $E$ be  a semistable vector 
bundle over $X$ of rank $n_0$ and degree $d_0$. The generalised Picard
bundle ${\mathcal W}_{\xi}(E)$ is by definition
the vector bundle over ${\mathcal M}_\xi$ defined by the direct
image $p_{\mathcal M_{\xi} *}({\mathcal U}_{\xi} \otimes p^*_X E)$
where
${\mathcal U}_{\xi}$
is a universal vector bundle over $X\times {\mathcal
M}_{\xi}$.
We obtain an inversion formula allowing us to recover $E$ from 
${\mathcal W}_{\xi}(E)$ and show that
the space of infinitesimal deformations of ${\mathcal
W}_{\xi}(E)$ is isomorphic to $H^1(X,\,End(E))$. 
This construction gives a locally complete family
of vector bundles over ${\mathcal M}_\xi$
parametrised by the moduli space ${\mathcal M}(n_0,d_0)$ of stable bundles
of rank $n_0$ and degree $d_0$  over $X$.
If $(n_0,d_0)=1$ and ${\mathcal W}_{\xi}(E)$ is stable for all 
$E\in{\mathcal M}(n_0,d_0)$, the construction determines an isomorphism
from ${\mathcal M}(n_0,d_0)$ to
a connected component ${\mathcal M}^0$ of a moduli space of stable sheaves 
over ${\mathcal M}_{\xi}$. This applies in particular when $n_0=1$, in which
case ${\mathcal M}^0$
is isomorphic to the Jacobian $J$ of $X$ as a polarised variety.
The paper as a whole is a generalisation of results of Kempf and Mukai on
Picard bundles over $J$.
\end{abstract}

\maketitle

\section{Introduction}

Let $X$ be a connected
nonsingular projective algebraic curve of genus $g \geq 2$ defined
over the complex numbers.
Let
$J$ denote the Jacobian (Picard variety) of $X$ and $J^d$ the variety
of line bundles of degree $d$ over $X$; thus in particular $J^0=J$.
Suppose $d\ge 2g-1$ and let $\mathcal L$ be a Poincar\'e (universal)
bundle over $X\times J^d$.  If we denote by $p_J$ the natural
projection from $X\times J^d$ to $J^d$, the direct image
$p_{J*}\mathcal L$ is then locally free and is called the {\em Picard
bundle} of degree $d$.

These bundles have been investigated by a number of authors over at
least the last 40 years. It may be noted that the projective bundle
corresponding to $p_{J*}\mathcal L$ can be identified with the
$d$-fold symmetric product $S^d(X)$. Picard bundles were studied in
this light by  A. Mattuck \cite{ma1, ma2} and I. G. Macdonald
\cite{mc} among others; both Mattuck and Macdonald gave formulae for
their Chern classes. Somewhat later R. C. Gunning \cite{gu1, gu2} gave
a more analytic treatment involving theta-functions. Later still, and
of especial relevance to us, G. Kempf \cite{ke} and S. Mukai
\cite{mu1} independently studied the deformations of the Picard
bundle; the problem then is to obtain an inversion formula showing
that all deformations of
$p_{J*}\mathcal L$ arise in a natural way. Kempf and Mukai proved
that $p_{J*}\mathcal L$ is simple and that, if $X$ is not
hyperelliptic, the space of infinitesimal deformations of 
$p_{J*}\mathcal L$ has dimension given by
$$
\dim\, H^1(J^d, \, End
(p_{J*}{\mathcal L}))=
2g.
$$
Moreover, all the infinitesimal deformations arise from genuine
deformations. In fact there is a complete family of deformations of
$p_{J*}\mathcal L$ parametrised by $J\times {\rm Pic}^0(J^d)$, the
two factors corresponding respectively to translations  in $J^d$ and
deformations of $\mathcal L$ (\cite[\S9]{ke}, \cite[Theorem
4.8]{mu1}). (The deformations of $\mathcal L$ are given by ${\mathcal
L}\mapsto{\mathcal L}\otimes p_J^* L$ for $L\in{\rm
Pic}^0(J^d)$.) Since $J$ is a principally polarised abelian variety
and $J^d\cong J$, ${\rm Pic}^0(J^d)$ can be identified with $J$
(strictly speaking ${\rm Pic}^0(J^d)$ is the dual abelian variety,
but the principal polarisation allows the identification).

Mukai's paper \cite{mu1} was set in a more general context involving a
transform which provides an equivalence between the derived category
of the category of ${\mathcal O}_A$-modules over an abelian variety
$A$ and the corresponding derived category on the dual abelian
variety $\hat A$. This technique has come to be known as the
Fourier--Mukai transform and has proved very useful in studying
moduli spaces of sheaves on abelian varieties and on some other
varieties.

Our object in this paper is to generalise the results of Kempf and
Mukai on deformations of Picard bundles to the moduli spaces of
higher rank vector bundles over $X$ with fixed determinant. In
particular we obtain an inversion formula for our generalised
Picard bundles and compute their spaces of infinitesimal deformations.
We also identify a family of
deformations which is locally complete and frequently globally complete
as well. The construction of the generalised Picard bundles together with 
the inversion formula
can be seen as a type of Fourier--Mukai transform.

We fix a holomorphic line bundle $\xi$ over $X$ of degree $d$.
Let ${\mathcal M}_{\xi}\, :=\, {\mathcal M}_{\xi}(n,d)$ be the
moduli space of stable vector bundles $F$ over $X$ with ${\rm
rank}(F) =n\ge2$, ${\rm deg}(F) =d$ and $\bigwedge^n F =\xi$. 
We assume that $n$ and $d$ are coprime, ensuring
the smoothness and completeness of ${\mathcal M}_{\xi}$, and that
$g\ge3$. We  assume also that if $g=3$ then  $n\geq 4$ and if $g=4$
then
$n\geq 3$. The case $g=2$ together with the three special cases $g=3$
with
$n=2,3$ and
$g=4$ with
$n=2$ are omitted in our main results since the method of proof
does not cover these cases.

It is known that there is a universal vector bundle over
$X\times {\mathcal M}_{\xi}$. Two such universal bundles
differ by tensoring with the pullback of a line bundle on
${\mathcal M}_{\xi}$. However, since ${\rm Pic}({\mathcal
M}_{\xi}) \, =\, {\mathbb Z}$, it is possible to choose canonically
a universal bundle. Let $l$ be the smallest positive number such
that $ld \equiv 1$ mod $n$. There is a unique universal vector
bundle ${\mathcal U}_{\xi}$ over $X\times {\mathcal M}_{\xi}$
such that ${\bigwedge}^n {\mathcal U}_{\xi}\big\vert_{\{x\}\times
{\mathcal M}_{\xi}} \,=\, {\Theta}^{\otimes l}$ \cite{Ra}, where
$x\in X$ and $\Theta$ is the ample generator of ${\rm
Pic}({\mathcal M}_{\xi})$.  Henceforth, by a universal bundle we
will always mean this canonical one. We denote by $p_X$ and
$p_{{\mathcal M}_{\xi}}$ the natural projections of $X\times{\mathcal
M}_{\xi}$ onto the two factors.

Now suppose that $n_0$ and $d_0$ are integers with $n_0\ge1$ and
\begin{equation}\label{deg}nd_0+n_0d>nn_0(2g-2).\end{equation}

For any semistable vector bundle $E$ of rank $n_0$ and degree $d_0$ over $X$,
let
$$
{\mathcal W}_{\xi}(E)\, :=\, p_{{\mathcal M}_{\xi}*} ({\mathcal
U}_{\xi}\otimes p^*_X E)
$$
be the direct image. The assumption \eqref{deg} ensures that
${\mathcal W}_{\xi}(E)$ is a locally free sheaf on
${\mathcal M}_{\xi}$ and all the higher direct images
of ${\mathcal U}_{\xi}\otimes p^*_X E$ vanish.
The rank of ${\mathcal W}_{\xi}(E)\,$ is $nd_0+n_0d+nn_0(1-g)$ and 
$H^i({\mathcal M}_{\xi},\, {\mathcal W}_{\xi}(E))\,\cong\, 
H^i(X\times {\mathcal M}_{\xi}, \,{\mathcal U}_{\xi}\otimes p^*_X
E)$. We shall refer to the bundles
${\mathcal W}_{\xi}(E)$ as {\em generalised Picard bundles}.

Our first main result is an inversion formula for this construction.

\medskip{\bf Theorem \ref{2.7}.}
{\em Suppose that \eqref{deg} holds and that $E$ 
is a semistable bundle of rank $n_0$ and degree $d_0$. Then 
$$E\cong{\mathcal R}^1p_{X*}(p_{{\mathcal M}_{\xi}}^*{\mathcal W}_{\xi}(E)
\otimes{\mathcal U}_{\xi}^*\otimes p_X^*K_X).$$}

\medskip
Following this, we show that, if $E$ is simple as well as semistable, then
${\mathcal W}_{\xi}(E)$ is simple (Corollary \ref{2.8}). Moreover

\medskip
{\bf Theorem \ref{2.9}.}
{\em Suppose that \eqref{deg} holds. 
For any semistable bundle $E$ of rank $n_0$ and degree $d_0$,
the space of infinitesimal deformations of the vector bundle
${\mathcal W}_{\xi}(E)$, namely $H^1({\mathcal M}_\xi ,\, End
( {\mathcal W}_{\xi}(E)))$, is canonically isomorphic to
$H^1(X,\, End E)$. In particular, if $E$ is also simple, 
$$
\dim H^1({\mathcal M}_\xi ,\, End  ( {\mathcal W}_{\xi}(E)))
\, =\, n_0^2(g-1)+1.
$$}

\medskip 

In the special case where $n=2$ and $E={\mathcal O}_X$, this was
proved by V. Balaji and P. A. Vishwanath in \cite{bv} using a
construction of M. Thaddeus \cite{ta}. For all $n$, it is known that
${\mathcal W}_{\xi}({\mathcal O}_X)$ is simple \cite{bef} and indeed
that it is stable (with respect to the unique polarisation of
${\mathcal M}_{\xi}$)
\cite{bbgn}; in fact the proof of stability generalises easily to
show that ${\mathcal W}_{\xi}(L)$ is stable for any line bundle $L$
for which \eqref{deg} holds. 
In this context, note
that Y. Li \cite{li} has proved a stability result for Picard
bundles over the non-fixed determinant moduli space ${\mathcal
M}(n,d)$, but this does not imply the result for ${\mathcal M}_{\xi}$.

If $(n_0,d_0)=1$, we can consider the bundles $\{ {\mathcal
W}_{\xi}(E)\}$ as a family of bundles over ${\mathcal M}_{\xi}$,
parametrized by
${\mathcal M}(n_0,d_0)$. We prove that this family is locally complete, i. e. that
the infinitesimal deformation map
$$
H^1(X,\, {\mathcal O}_X)\, \longrightarrow\, H^1({\mathcal M}_{\xi},\,
End({\mathcal W}_{\xi}(E)))
$$
is an isomorphism for all $E$ (Theorem \ref{5.1}). If all
the bundles ${\mathcal W}_{\xi}(E)$ are stable, the family is also
globally complete (Theorem \ref{5.3}). Finally, in the case $n_0=1$, we
obtain an isomorphism of polarised varieties between $J$ and a connected
component of a moduli space of bundles over ${\mathcal M}_{\xi}$ 
(Theorem \ref{5.4}), which in turn leads to a Torelli theorem 
(Corollary \ref{5.5}).

The layout of the paper is as follows. In sections 2, 3 and 4, we obtain 
cohomological results. The techniques are quite similar to those of Kempf 
and Mukai except for our use of Hecke transformations; however our
calculations are more complicated since we cannot
exploit the special properties of abelian varieties. In sections 5, 6 and 7,
we then use these results to obtain our main theorems.

{\it Notation and assumptions.} We work throughout over the
complex numbers and suppose that $X$ is a connected nonsingular
projective algebraic curve of genus $g\ge3$. We suppose that $n\ge2$
and that if
$g=3$ then $n\ge4$ and if $g=4$ then
$n\ge3$. We assume moreover that $(n,d)=1$ and that \eqref{deg} holds. In
general, we denote the natural projections of
$X\times Y$ onto its factors by $p_X$, $p_Y$.
For a variety $X\times Y\times
Z$, we denote by $p_i\ (i=1,2,3)$ the projection onto the $i$-th
factor and by $p_{ij}$ the projection onto the Cartesian
product of the $i$-th and the
$j$-th factors. Finally, for any $x\in X$, we denote by ${\mathcal
U}_x$ the bundle over ${\mathcal M}_{\xi}$ obtained by restricting
${\mathcal U}_{\xi}$ to $\{x\}\times{\mathcal M}_{\xi}$.

\smallskip
{\it Acknowledgements.} Part of this work was done 
during a visit of the authors to the Abdus Salam ICTP, Italy. The 
authors thank ICTP for its hospitality. The third author would like
also to thank the Isaac Newton Institute, Cambridge and the organisers
of the HDG programme for their hospitality during the preparation of
an earlier version of this paper.
\smallskip

\section{Cohomology of ${\mathcal W}_{\xi}(E_1)\otimes{\mathcal
W}_\xi(E_2)^*$}

Our principal object in this section and the two following sections 
is to compute the cohomology groups $H^i({\mathcal M}_{\xi},
{\mathcal W}_{\xi}(E_1)\otimes{\mathcal W}_\xi(E_2)^*)$ for $i=0,1$, where
$E_1$ and $E_2$ are semistable bundles of rank $n_0$ and degree $d_0$ 
satisfying \eqref{deg}.

\begin{propn}\label{2.1}
Suppose that \eqref{deg} holds and $E_1$ and $E_2$ are semistable bundles of 
rank $n_0$ and degree $d_0$. Then
\begin{eqnarray*}
H^i({\mathcal M}_\xi,
\, {\mathcal W}_{\xi}(E_1)\otimes {\mathcal W}_{\xi}(E_2)^* )& \cong &
H^i(X\times {\mathcal M}_\xi,\, 
{\mathcal U}_{\xi}\otimes p_X^*E_1\otimes p_{{\mathcal
M}_{\xi}}^*{\mathcal W}_\xi (E_2)^*)\\
&\cong& H^{i+1}(X\times {\mathcal M}_\xi \times X, \, p_{12}^*
{\mathcal
U}_{\xi}\otimes p^*_1E_1 \otimes p^*_{23}{\mathcal U}^*_\xi\otimes 
p^*_3 E_2^*\otimes  p^*_3K_X)
\end{eqnarray*}
for $i \geq 0$, where $K_X$ is the canonical line bundle over $X$.
\end{propn}

\begin{proof}
Recall first that, if $E$ and $F$ are semistable, then so is
$F\otimes E$. It then follows from \eqref{deg} that  $H^1(X,\, 
{\mathcal U}_{\xi}\big\vert_{X\times \{v\}}\otimes E_1) =0$ for all $v\in
{\mathcal M}_\xi$.  Using the projection formula and the Leray
spectral sequence we have
$$
H^i(X\times {\mathcal M}_\xi,\, 
{\mathcal U}_{\xi}\otimes p^*_X E_1\otimes p_{{\mathcal
M}_{\xi}}^*{\mathcal W}_\xi (E_2)^*)\, \cong\, H^i({\mathcal M}_\xi,
\, {\mathcal W}_{\xi}(E_1)\otimes {\mathcal W}_{\xi}(E_2)^* )\, .
$$
This proves the first isomorphism.

In the same way, \eqref{deg} gives
$$
H^0(X,\, ({\mathcal U}_{\xi}\big\vert_{X\times \{v\}})^*\otimes
E^*_2\otimes K_X) \cong
H^1(X,\, {\mathcal U}_{\xi}\big\vert_{X\times \{v\}}\otimes E_2)^* \, =\, 0\, .
$$
Consequently, the projection formula gives
$$
{\mathcal R}^{i} p_{12*}(p_{23}^*({\mathcal U}^*_\xi)\otimes
p_3^*E^*_2 \otimes p_3^* K_X) \, =\, 0
$$
for $i\not=1$, and we have by relative Serre duality
$$
{\mathcal R} ^{1} p_{12*}(p_{23}^*({\mathcal U}^*_\xi)\otimes
p_3^*E^*_2 \otimes p_3^*K_X)\, \cong\, p_{{\mathcal M}_{\xi}}^* {\mathcal
W}_\xi (E_2)^*\, .
$$
Finally, using the projection formula and the Leray spectral
sequence, it follows that
$$
H^i(X\times {\mathcal M}_\xi,\,
{\mathcal U}_{\xi}\otimes p^*_XE_1 \otimes 
p_{{\mathcal M}_{\xi}}^*{\mathcal W}_\xi (E_2)^*)
$$
$$
\cong\,  
H^{i+1}(X\times {\mathcal M}_\xi \times X, \, p_{12}^*
{\mathcal U}_{\xi} \otimes p^*_1E_1\otimes
  p^*_{23}{\mathcal U}^*_\xi \otimes p^*_3E^*_2 \otimes p^*_3 K_X).
$$
for $i\geq 0$. This completes the proof.
\end{proof}

\begin{rem}\label{2.2}
\begin{em}
Proposition \ref{2.1} can be formulated
in a more general context.
Let $V_1$, $V_2$ be flat families of vector
bundles over $X$
parametrised by a complete irreducible variety $Y$ 
such that for each $y\in Y$ we have
$H^1(X,\, V_i|_{X\times \{y\}})\,=\,0$ for $i=1,2$. Under this
assumption
$$
H^i(Y,\, p_{Y*}V_1\otimes (p_{Y*}V_2)^*)\,\cong\, H^{i+1}(X\times Y
\times X,
\, p^*_{12}V_1 \otimes p^*_{23}V_2^*\otimes p^*_3K_X)\, .
$$
The proof is the same as for Proposition \ref{2.1}.
\end{em}\end{rem}

We can now state the key result which enables us to calculate the
cohomology groups in which we are interested. We
denote by ${\mathcal R}^i$ the $i$-th direct image of $ p^*_{12}
{\mathcal U}_{\xi}\otimes p^*_1E_1 \otimes
p^*_{23}{\mathcal U}^*_\xi$ for the projection $p_3$, that is 
$${\mathcal R}^i := {\mathcal R}^i p_{3*}( p^*_{12}
{\mathcal U}_{\xi}\otimes p^*_1E_1\otimes
p^*_{23}{\mathcal U}^*_\xi).$$

\begin{propn}\label{leray} For $i\ge0$, there exists an exact sequence
$$
0\,\longrightarrow\, H^1(X,\,{\mathcal R}^i\otimes E^*_2\otimes K_X)
\,\longrightarrow\,
H^{i}({\mathcal M}_\xi , {\mathcal W}_{\xi}(E_1)\otimes
{\mathcal W}_\xi(E_2)^*)
$$
\begin{equation}\label{1}
\hspace{5.5cm}\, \longrightarrow\, H^0(X,{\mathcal R}^{i+1}\otimes E^*_2\otimes
K_X)\,\longrightarrow \,0.
\end{equation}
\end{propn}
\begin{proof} Since $\dim X=1$, the Leray spectral sequence for $p_3$ gives
$$
0\, \longrightarrow\, H^1(X,{\mathcal R}^i\otimes E^*_2\otimes K_X)
\, \longrightarrow\,
H^{i+1}(X\times {\mathcal M}_\xi \times X, \, p^*_{12}
{\mathcal U}_{\xi}\otimes p^*_1E_1 \otimes
p^*_{23}{\mathcal U}^*_\xi\otimes p^*_3 E^*_2\otimes p^*_3K_X)
$$
\begin{equation}\label{2}
\hspace{5.5cm}\, \longrightarrow\, H^0(X,{\mathcal R}^{i+1}\otimes E^*_2
\otimes K_X)\, \longrightarrow\, 0\, .
\end{equation}
The result now follows at once from Proposition \ref{2.1}.
\end{proof}

In order to use this proposition, we must compute the ${\mathcal R}^i$
for $i\le2$. We can already prove

\begin{propn}\label{2.3}
${\mathcal R}^0 \,=\, 0$.
\end{propn}

\begin{proof}
Note that for any $x \in X$
\begin{equation}\label{3}
H^0(X\times {\mathcal M}_{\xi},\, {\mathcal U}_{\xi} \otimes
p^*_XE_1
\otimes p_{{\mathcal M}_{\xi}}^*{\mathcal U}^*_{x}) \cong H^0(X ,\,
E_1\otimes  p_{X*}({\mathcal U}_{\xi} 
\otimes p_{{\mathcal M}_{\xi}}^*{\mathcal U}^*_{x}))\, .
\end{equation}

{}From \cite{ib} and \cite{nr} we know that for generic $y \in X$, 
the two vector bundles ${\mathcal U}_{y}$
and ${\mathcal U}_{x}$ are non-isomorphic and stable. Hence
$H^0({\mathcal M}_{\xi}, \,{\mathcal U}_{y} \otimes {\mathcal
U}^*_{x})=0.$ This implies that
\begin{equation}\label{4}
p_{X*}({\mathcal U}_{\xi} 
\otimes p_{{\mathcal M}_{\xi}}^*{\mathcal U}^*_{x}) =0\, .
\end{equation}
So \eqref{3} gives
$$
{\mathcal R}^0:= p_{3*}(p_{12}^* {\mathcal
U}_{\xi} \otimes p^*_1E_1\otimes p^*_{23}{\mathcal U}^*_\xi) \,=\,
0
$$
and the proof is complete.
\end{proof}

In the next two sections we will use Hecke transformations and
a diagonal argument to show that $R^2=0$ and to compute $R^1$.

\section{The Hecke Transformation}

In this section we will use Hecke transformations to compute 
the cohomology groups $H^i(X\times{\mathcal M}_\xi ,\,
{\mathcal U}_{\xi}\otimes p^*_1E_1\otimes 
p_{{\mathcal M}_{\xi}}^*{\mathcal U}^*_x)$ for any $x\in X$.
The details of the Hecke transformation and its properties
can be found in \cite{nr,nr2}. 
We will briefly describe it and
note those properties that will be needed here.

Fix a point $x\in X$. Let ${\mathbb P}({\mathcal U}_x)$ denote the
projective bundle over ${\mathcal M}_\xi$ consisting of lines in
${\mathcal U}_x$. If $f$ denotes the natural projection of
${\mathbb P}({\mathcal U}_x)$ to ${\mathcal M}_\xi$ and 
${\mathcal O}_{{\mathbb P}({\mathcal U}_x)}(-1)$  the
tautological line bundle then 
$$f_*{\mathcal O}_{{\mathbb P}({\mathcal U}_x)}(1) \, 
\cong\, {\mathcal U}^*_x\, ,
$$
and ${\mathcal R}^j f_* {\mathcal O}_{{\mathbb P}
({\mathcal U}_x)}(1) \, = 0$
for all $j>0.$ From the commutative diagram
$$
\begin{array}{ccc}
X\times {\mathbb P}({\mathcal U}_x)&
\stackrel{{\rm Id}_X\times f}{\longrightarrow}&X \times {\mathcal M}_\xi\\
p_{{\mathbb P}({\mathcal U}_x)}\Big\downarrow&&\Big\downarrow p_{\mathcal
M}\\ {\mathbb P}({\mathcal U}_x)&
\stackrel{f}{\longrightarrow}&{\mathcal M}_\xi
\end{array}
$$
and the base change theorem, we
 deduce that
$$H^i(X\times{\mathcal M}_\xi ,\,
{\mathcal U}_{\xi}\otimes p^*_XE_1\otimes p_{{\mathcal
M}_{\xi}}^*{\mathcal U}^*_x)$$
\begin{equation}\label{6}
\cong
H^i(X\times {\mathbb P}({\mathcal U}_x),\, ({\rm Id}_X\times f)^*
({\mathcal U}_{\xi}\otimes p^*_XE_1)
\otimes p_{{\mathbb P}({\mathcal U}_x)}^*
{\mathcal O}_{{\mathbb P}({\mathcal
U}_x)}(1))
\end{equation}
for all $i$. 

Moreover, since $p_{{\mathbb P}({\mathcal U}_x)*}({\rm Id}_X\times f)^*
({\mathcal U}_{\xi}\otimes p^*_XE_1) \, \cong\, f^*{\mathcal
W}_\xi(E_1)$, there is a canonical isomorphism
$$
H^i(X\times {\mathbb P}({\mathcal U}_x),\, ({\rm Id}_X\times f)^*
({\mathcal U}_{\xi}\otimes p^*_XE_1 )
\otimes p_{{\mathbb P}({\mathcal U}_x)}^*{\mathcal O}_{{\mathbb
P}({\mathcal U}_x)}(1))$$
\begin{equation}\label{7}
\cong\,
H^i({\mathbb P}({\mathcal U}_x),\, f^*{\mathcal
W}_\xi(E_1)\otimes {\mathcal O}_{{\mathbb P}({\mathcal U}_x)}(1))
\end{equation}
for all $i$.

To compute the cohomology groups 
$H^i({\mathbb P}({\mathcal U}_x),\, f^*{\mathcal
W}_\xi(E_1)\otimes {\mathcal O}_{{\mathbb P}({\mathcal U}_x)}(1))$ 
we use Hecke transformations.

A point in ${\mathbb P}({\mathcal U}_x)$ 
represents a stable vector bundle
$F$ and a line $l$ in the fibre $F_x$ at $x$, or equivalently a
non-trivial exact sequence
\begin{equation}\label{8}
0\, \longrightarrow\, F\, \longrightarrow\, F' \,
\longrightarrow\, {\mathbb C}_x \, \longrightarrow\, 0 \, 
\end{equation}
determined up to a scalar multiple; here ${\mathbb C}_x$ denotes the
torsion sheaf supported at $x$ with stalk $\mathbb C$. The sequences
\eqref{8} fit together to form a universal sequence
\begin{equation}\label{9}
0\, \longrightarrow\, ({\rm Id}_X\times f)^*
{\mathcal U}_{\xi}\otimes p_{{\mathbb P}({\mathcal U}_x)}^*{\mathcal
O}_{{\mathbb P}({\mathcal U}_x)}(1)\, \longrightarrow\, {\mathcal F} \,
\longrightarrow\, p_X^*{\mathbb C}_x \, \longrightarrow\, 0 \, 
\end{equation}
on $X\times{\mathbb P}({\mathcal U}_x)$.
  If $\eta$ denotes the line bundle $\xi\otimes {\mathcal O}_X(
x)$ over $X$ and ${\mathcal
M}_\eta $ the moduli space of stable bundles $ {\mathcal M}_\eta(n,d+1)$ 
then from \eqref{8} and \eqref{9} we get a rational map
$$
\gamma\, :\, {\mathbb P}({\mathcal U}_x)\, - - \rightarrow\,{\mathcal
M}_\eta$$
which sends any pair $(F,l)$ to $F'$. This map is not everywhere
defined since the bundle $F'$ in \eqref{8} need not be stable. 

Our next
object is to find a Zariski-open subset
$Z$ of ${\mathcal M}_\eta$, over which $\gamma$ is defined and is a
projective fibration, such that the complement of $Z$ in ${\mathcal
M}_\eta$ has codimension at least 4. The construction and
calculations are similar to those of \cite[Proposition 6.8]{nr}, but
our results do not seem to follow directly from that proposition.

As in \cite[\S8]{nr} or \cite[\S5]{nr2}, we define a bundle $F'$ to
be $(0,1)$-{\it stable} if, for every proper subbundle $G$ of $F'$,
$$\frac{\deg G}{\hbox{rk}G}<\frac{\deg F'-1}{\hbox{rk}F'}.$$
Clearly every $(0,1)$-stable bundle is stable. We denote by $Z$ the
subset of ${\mathcal M}_{\eta}$ consisting of $(0,1)$-stable bundles. 

\begin{lemma}\label{3.1}
{\rm (i)} $Z$ is a Zariski-open subset of ${\mathcal
M}_{\eta}$ whose complement has codimension at least $4$.

{\rm (ii)} $\gamma$ is a projective fibration over $Z$ and
$\gamma^{-1}(Z)$ is a Zariski-open subset of ${\mathbb P}({\mathcal
U}_x)$ whose complement has codimension at least $4$.
\end{lemma}

\begin{proof}
(i) The fact that $Z$ is Zariski-open is standard (see
\cite[Proposition 5.3]{nr2}).

The bundle $F'\in{\mathcal M}_{\eta}$ of rank $n$ and degree $d+1$
fails to be $(0,1)$-stable if and only if it has a subbundle $G$ of
rank $r$ and degree $e$ such that $ne\ge r((d+1)-1)$, i.e.,
\begin{equation}\label{10}
rd\,\le\, ne\, .
\end{equation}
By considering the extensions
$$0\longrightarrow G\longrightarrow F'\longrightarrow H\longrightarrow
0,$$
we can estimate the codimension of ${\mathcal M}_{\eta}-Z$ and show
that it is at least
\begin{equation}\label{11}
\delta=r(n-r)(g-1)+(ne-r(d+1))
\end{equation}
(compare the proof of \cite[Proposition 5.4]{nr2}). Note that, since
$(n,d)=1$, \eqref{10} implies that $rd\le ne-1$. Given that $g\ge3$, we
see that $\delta<4$ only if $g=3$, $n=2,3$ or $g=4$, $n=2$. These are
exactly the cases that were excluded in the introduction.

(ii) $\gamma^{-1}(Z)$ consists of all pairs $(F,l)$ for which the
bundle $F'$ in \eqref{8} is $(0,1)$-stable. As in (i), this is a
Zariski-open subset.
It follows at once from \eqref{10} that, if $F'$ is $(0,1)$-stable,
then $F$ is stable. So, if $F'\in Z$, it follows from
\eqref{8} that $\gamma^{-1}(F')$ can be identified with the projective
space ${\mathbb P}(F'_x{}^*)$. Using the universal projective bundle on
$X\times{\mathcal M}_{\eta}$, we see that $\gamma^{-1}(Z)$ is a
projective fibration over $Z$ (not necessarily locally trivial).

Suppose now that $(F,l)$ belongs to the complement of
$\gamma^{-1}(Z)$ in ${\mathbb P}({\mathcal U}_x)$. This means that the
bundle $F'$ in \eqref{8} is not $(0,1)$-stable and therefore possesses a
subbundle $G$ satisfying \eqref{10}. If $G\subset F$, this contradicts the
stability of $F$. So there exists an exact sequence
$$0\longrightarrow G'\longrightarrow G\longrightarrow{\mathbb
C}_x\longrightarrow 0$$
with $G'$ a subbundle of $F$ of rank $r$ and degree $e-1$. Moreover,
since $G$ is a subbundle of $F$, $G'_x$ must contain the line $l$.
For fixed $r,e$, these conditions determine a subvariety of ${\mathbb
P}({\mathcal U}_x)$ of dimension at most
$$ 
(r^2(g-1)+1)+((n-r)^2(g-1)+1)-g+(r-1)
+((g-1)r(n-r)+(rd-n(e-1))-1).
$$
Since $\dim{\mathbb
P}({\mathcal U}_x)=n^2(g-1)-g+n$, a simple calculation shows that the
codimension is at least the number $\delta$ given by \eqref{11}. As in
(i),
this gives the required result.
\end{proof}

By Lemma \ref{3.1}(ii) and a Hartogs-type theorem (see \cite[Theorem
3.8 and Proposition 1.11]{hart}) we have an isomorphism
\begin{equation}\label{12}
H^i({\mathbb P}({\mathcal U}_x),\, f^*{\mathcal
W}_\xi(E_1)\otimes {\mathcal O}_{{\mathbb P}({\mathcal U}_x)}(1))
\cong H^i(\gamma ^{-1}(Z),\, f^*{\mathcal
W}_\xi(E_1)\otimes {\mathcal O}_{{\mathbb P}
({\mathcal U}_x)}(1)|_{\gamma^{-1}(Z)})
\end{equation}
for $i\le2$.

Now let $F'\in Z$. As in the proof of Lemma \ref{3.1}, we identify
$\gamma^{-1}(F')$ with ${\mathbb P}(F'_x{}^*)$ and denote it by ${\mathbb P}$.
On $X\times{\mathbb P}$ there is a universal exact sequence
\begin{equation}\label{13}
0\longrightarrow{\mathcal E}\longrightarrow
p_X^*F'\longrightarrow p_{\mathbb P}^*{\mathcal O}_{\mathbb P}(1)\otimes
p_X^*{\mathbb C}_x\longrightarrow 0\, .
\end{equation}
The
restriction of \eqref{13} to any point of ${\mathbb P}$ is isomorphic to
the
corresponding sequence \eqref{8}.

\begin{propn}\label{3.2}
Let ${\mathcal F}$ be defined by the universal sequence
\eqref{9}. Then
$${\mathcal F}\big\vert_{X\times{\mathbb P}}\,\cong\, p_X^*F'\otimes
p_{\mathbb
P}^*{\mathcal O}_{\mathbb P}(-1)\, .$$
\end{propn}

\begin{proof}
Restricting \eqref{9} to $X\times{\mathbb P}$
gives
$$
0\, \longrightarrow\, ({\rm Id}_X\times f)^*
{\mathcal U}_{\xi}\otimes p_{{\mathbb P}({\mathcal U}_x)}^*{\mathcal
O}_{{\mathbb P}({\mathcal U}_x)}(1))\big\vert_{X\times{\mathbb P}}
\longrightarrow\, {\mathcal F}\big\vert_{X\times{\mathbb P}}
\,
\longrightarrow\, p_X^*{\mathbb C}_x \, \longrightarrow\, 0.
$$
This must coincide with the universal sequence
\eqref{13} up to tensoring
by some line bundle lifted from ${\mathbb P}$. The result follows.
\end{proof}

Next we tensor \eqref{9} by $p_X^*E_1$, restrict it to
$X\times\gamma^{-1}(Z)$ and take the direct image on
$\gamma^{-1}(Z)$. This gives
\begin{equation}\label{14}
0\, \longrightarrow\, f^*
{\mathcal W}_{\xi}(E_1)\otimes {\mathcal O}_{{\mathbb P}({\mathcal
U}_x)}(1)\big\vert_{\gamma^{-1}(Z)} \longrightarrow\, p_{{\mathbb
P}({\mathcal U}_x)*}({\mathcal F}\otimes
p_X^*E_1)\big\vert_{\gamma^{-1}(Z)}
\,
\longrightarrow\, {\mathcal O}_{\gamma^{-1}(Z)}^{\oplus n_0} \, 
\longrightarrow\,0\, .
\end{equation}

\begin{propn}\label{3.3}
$R^i_{\gamma*}(p_{{\mathbb P}({\mathcal
U}_x)*}({\mathcal F}\otimes p_X^*E_1)\big\vert_{\gamma^{-1}(Z)})=0$
for all $i$.
\end{propn}

\begin{proof}
It is sufficient to show that $p_{{\mathbb P}({\mathcal
U}_x)*}({\mathcal F}\otimes p_X^*E_1)\big\vert_{\mathbb P}$ has trivial
cohomology. By Proposition \ref{3.2}
\begin{eqnarray*}p_{{\mathbb P}({\mathcal U}_x)*}({\mathcal F}\otimes
p_X^*E_1)\big\vert_{\mathbb P}&\cong& p_{{\mathbb
P}*}(p_X^*F'\otimes p_X^*E_1\otimes p_{\mathbb P}^*{\mathcal
O}_{\mathbb P}(-1))\\&\cong& H^0(X,F'\otimes
E_1)\otimes{\mathcal O}_{\mathbb P}(-1)\end{eqnarray*}
and the result follows.
\end{proof}

\begin{coro}\label{3.4}
${\mathcal R}^i_{\gamma_*} (f^*{\mathcal
W}_\xi(E_1)\otimes {\mathcal O}_{{\mathbb P}({\mathcal
U}_x)}(1)\big\vert_{\gamma^{-1}(Z)})
\, = 0$ for $i \ne1$. Moreover, 
$${\mathcal R}^1_{\gamma_*} (f^*{\mathcal
W}_\xi(E_1)\otimes {\mathcal O}_{{\mathbb P}({\mathcal
U}_x)}(1)\big\vert_{\gamma^{-1}(Z)})
\cong{\mathcal O}_Z^{\oplus n_0}.$$
\end{coro}

\begin{proof}
This follows at once from \eqref{14} and Proposition \ref{3.3}.
\end{proof}

Now we are in a position to compute the cohomology groups of
$H^i(X\times{\mathcal M}_\xi ,\,
{\mathcal U}_{\xi}\otimes p^*_XE_1\otimes 
p_{{\mathcal M}_{\xi}}^*{\mathcal U}^*_x)$ for  $i=1,2$.

\begin{propn}\label{3.5}
$H^2(X\times{\mathcal M}_\xi ,\,
{\mathcal U}_{\xi}\otimes p^*_XE_1\otimes 
p_{{\mathcal M}_{\xi}}^*{\mathcal U}^*_x)=0 $ for any $x\in X.$
\end{propn}

\begin{proof}
The combination of \eqref{6}, \eqref{7} and \eqref{12}
yields
$$
H^2(X\times{\mathcal M}_\xi ,\,
{\mathcal U}_{\xi}\otimes p^*_XE_1\otimes 
p_{{\mathcal M}_{\xi}}^*{\mathcal U}^*_x) \,\cong\, 
H^2(\gamma ^{-1}(Z),\, f^*{\mathcal
W}_\xi(E_1)\otimes {\mathcal O}_{{\mathbb P}
({\mathcal U}_x)}(1)\big\vert_{\gamma^{-1}(Z)})\, .
$$
Using Corollary \ref{3.4} and Lemma \ref{3.1}(i), the Leray spectral
sequence
for the map $\gamma$ gives$$
\begin{array}{lll}
H^2(\gamma^{-1}(Z),\,   f^*{\mathcal
W}_\xi(E_1)\otimes {\mathcal O}_{{\mathbb P}
({\mathcal U}_x)}(1)\big\vert _{\gamma^{-1}(Z)})\,
&\cong\,& H^{1}(Z,\, {\mathcal O}_Z)^{\oplus n_0}\\
&\cong\,& H^{1}({\mathcal M}_\eta,\, {\mathcal O}_{{\mathcal
M}_{\eta}})^{\oplus n_0}.
\end{array}
$$
It is known that
$H^1({\mathcal M}_\eta,\, {\mathcal O}_{{\mathcal M}_{\eta}})\,=\,0$
\cite{dn}. Therefore, 
$$
H^2(X\times{\mathcal M}_\xi ,\,
{\mathcal U}_{\xi}\otimes p^*_XE_1\otimes 
p_{{\mathcal M}_{\xi}}^*{\mathcal U}^*_x) = 0\, .
$$
\end{proof}

We can now prove
\begin{propn}\label{2.4}
${\mathcal R}^2 \,=\,0.$  
\end{propn}

\begin{proof} This is
an immediate consequence of Proposition \ref{3.5}.
\end{proof}
\medskip

\begin{propn}\label{3.6}
For any point $x\in X$, $\dim H^1(X\times{\mathcal M}_\xi ,\,
{\mathcal U}_{\xi}\otimes p^*_XE_1\otimes 
p_{{\mathcal M}_{\xi}}^*{\mathcal U}^*_x)\, =\, n_0$.
\end{propn}

\begin{proof}
As in the proof of Proposition \ref{3.5} 
we conclude that
$$
H^1(X\times{\mathcal M}_\xi ,\,
{\mathcal U}_{\xi}\otimes p^*_XE_1\otimes 
p_{{\mathcal M}_{\xi}}^*{\mathcal U}^*_x)
\, \cong\,H^{0}({\mathcal M}_\eta,\, {\mathcal O}_{{\mathcal
M}_{\eta}})^{\oplus n_0}\, .
$$
Now ${\mathcal M}_{\eta}$ is just the non-singular part of the moduli
space of semistable bundles of rank $n$ and determinant $\eta$, and
the latter space is complete and normal. So
$\dim H^0({\mathcal M}_{\eta},\, {\mathcal O}_{{\mathcal
M}_{\eta}})\,=\, 1$.
\end{proof}

\begin{coro}\label{r1}
${\mathcal R}^1$ is a vector bundle of rank $n_0$.
\end{coro}

The complete identification of ${\mathcal R}^1$ will be given in the
next section.

\medskip
\begin{rem}\label{3.7}
\begin{em}
Since the fibres of $\gamma$ are projective spaces, we have
$\gamma_*{\mathcal O}_{\gamma^{-1}(Z)}\cong{\mathcal O}_Z$ and all the
higher direct images of ${\mathcal O}_{\gamma^{-1}(Z)}$ are $0$. Hence
$$H^i(Z,{\mathcal O}_Z)\cong H^i(\gamma^{-1}(Z),{\mathcal
O}_{\gamma^{-1}(Z)})$$for all $i$. Similarly
$$H^i({\mathbb P}({\mathcal U}_x),{\mathcal O}_{{\mathbb P}({\mathcal
U}_x)})\cong H^i({\mathcal M}_{\xi},{\mathcal O}_{{\mathcal
M}_{\xi}})=0$$ for
$i>0$ since
${\mathcal M}_{\xi}$ is a smooth projective rational variety. It
follows from the proof of Lemma \ref{3.1} that,
if we define $\delta$ as in
\eqref{11},
$$\delta\ge i+2\ge3\Rightarrow H^i(Z,{\mathcal O}_Z)=0.$$
The proof of Proposition \ref{3.5} now gives
\begin{equation}\label{del}
\delta\ge i+2\ge4\Rightarrow H^i(X\times{\mathcal M}_{\xi},
{\mathcal U}_{\xi}\otimes p^*_XE_1\otimes p^*_{{\mathcal M}_{\xi}}{\mathcal
U}_x^*)=0.
\end{equation}
\end{em}\end{rem}
\begin{propn}\label{delta} Suppose that \eqref{deg} holds and that $E_1$, 
$E_2$ are semistable bundles of rank $n_0$ and degree $d_0$.
If $\delta\ge i+3\ge5$, then
$H^i({\mathcal M}_{\xi},{\mathcal
W}_{\xi}(E_1)\otimes {\mathcal W}_{\xi}(E_2)^*)=0.$
\end{propn}
\begin{proof} It follows from \eqref{del} that $R^i=R^{i+1}=0$. The result
now follows from Proposition \ref{leray}.\end{proof}

\begin{coro}\label{h2} Suppose that \eqref{deg} holds and that $E$ 
is a semistable bundle of rank $n_0$ and degree $d_0$. 
Then
$$H^2({\mathcal M}_{\xi},End({\mathcal W}_{\xi}(E)))=0$$
except possibly when $g=3, n=2, 3, 4$; $g=4,n=2$; $g=5,n=2$.
\end{coro}
\begin{proof} Take $E_1=E_2=E$ and $i=2$ in Proposition \ref{delta}.
We need to show that $\delta\ge5$. In fact
it follows from \eqref{12} that the exceptional cases are precisely
those for which $\delta<5$.\end{proof}

\section{A Diagonal Argument}

Let $\Delta $ be the diagonal divisor in $X\times X$. Pull back the 
exact sequence
$$0\, \longrightarrow\, {\mathcal O}(-{\Delta}) \, \longrightarrow\, 
{\mathcal O}\, \longrightarrow\, {\mathcal O}_{\Delta}\,
\longrightarrow\, 0
$$
to $X\times {\mathcal M}_{\xi} \times X$ and tensor it with 
$ p^*_{12}
{\mathcal U}_{\xi}\otimes p^*_1E_1 \otimes
p^*_{23}{\mathcal U}^*_\xi$. Now, the direct image 
sequence for the projection 
$p_3$ gives the following exact sequence over $X$
$$
\longrightarrow\, {\mathcal R}^i p_{3*} ( p^*_{12}
{\mathcal U}_{\xi}\otimes p^*_1E_1 \otimes
p^*_{23}{\mathcal U}^*_\xi (-{\Delta}))\, \longrightarrow\,
{\mathcal R}^i
 \, \longrightarrow\, {\mathcal R}^i p_{3*} ( p^*_{12}
{\mathcal U}_{\xi}\otimes p^*_1E_1 \otimes
p^*_{23}{\mathcal U}^*_\xi |_{\Delta \times {\mathcal M}_{\xi}})
$$
\begin{equation}\label{15}
\longrightarrow\, {\mathcal R}^{i+1} p_{3*} ( p^*_{12}
{\mathcal U}_{\xi}\otimes p^*_1E_1 \otimes
p^*_{23}{\mathcal U}^*_\xi (-{\Delta}))
\,\longrightarrow\,\ldots 
\end{equation}

The following proposition will be used in computing  
${\mathcal R}^1.$ 

\begin{propn}\label{4.1}
For any $E_1$, the direct images of
$$
p^*_{12}
{\mathcal U}_{\xi}\otimes p^*_1E_1 \otimes
p^*_{23}{\mathcal U}^*_\xi |_{\Delta \times {\mathcal M}_{\xi}}
$$
have the following description:

\begin{enumerate}
\item ${\mathcal R}^0 p_{3*}( p^*_{12}
{\mathcal U}_{\xi}\otimes p^*_1E_1 \otimes
p^*_{23}{\mathcal U}^*_\xi |_{\Delta \times {\mathcal M}_{\xi}})
\,\cong\, E_1$
\item ${\mathcal R}^1 p_{3*}( p^*_{12}
{\mathcal U}_{\xi}\otimes p^*_1E_1 \otimes
p^*_{23}{\mathcal U}^*_\xi |_{\Delta \times {\mathcal M}_{\xi}})
\,\cong\, E_1\otimes TX$
\item  ${\mathcal R}^2 p_{3*}( p^*_{12}
{\mathcal U}_{\xi}\otimes p^*_1E_1 \otimes
p^*_{23}{\mathcal U}^*_\xi |_{\Delta \times {\mathcal M}_{\xi}})\,=\,0$
\end{enumerate}
where $TX$ is the tangent bundle of $X$.
\end{propn}

\begin{proof}
Identifying $\Delta $ with $X$ we have
$${\mathcal R}^i p_{3*} ( p^*_{12}
{\mathcal U}_{\xi}\otimes p^*_1E_1 \otimes
p^*_{23}{\mathcal U}^*_\xi |_{\Delta \times {\mathcal M}_{\xi}})
\,\cong \, {\mathcal R}^i p_{X*}( 
{\mathcal U}_{\xi}\otimes p^*_XE_1 \otimes 
{\mathcal U}^*_\xi ).$$ 

The proposition follows from a result of 
Narasimhan and Ramanan \cite[Theorem 2]{nr} that says
\begin{equation}\label{16}
H^i({\mathcal M}_{\xi},\, {\mathcal U}_x \otimes {\mathcal U}^*_x )\,
\cong \,
\left\{\begin{array}{ll} {\mathbb C} &\mbox{if} \ \ i=0,1\\
0 & \mbox{if} \ \
i=2.  \end{array}\right.
\end{equation}
For $i=0$ the isomorphism is given
by the obvious inclusion of ${\mathcal O}_{{\mathcal M}_{\xi}}$ in 
${\mathcal U}_x
\otimes {\mathcal U}^*_x $ and therefore globalises to give
${\mathcal R}^0 p_{X*}(  {\mathcal U}_{\xi} \otimes 
{\mathcal U}^*_\xi )\cong{\mathcal O}_X$. Similarly for $i=1$ the
isomorphism is given by the infinitesimal deformation map of
${\mathcal U}_{\xi}$ regarded as a family of bundles over $\mathcal
M_{\xi}$ parametrised by $X$; this globalises to ${\mathcal R}^1
p_{X*}({\mathcal U}_{\xi}\otimes 
{\mathcal U}^*_\xi )\cong TX$.
\end{proof}

Propositions \ref{4.1} and \ref{2.3} and the exact sequence \eqref{15}
together give the following exact sequence of direct images
$$
0\,\longrightarrow\, E_1 \,\longrightarrow\, {\mathcal R}^1 p_{3*}(
p^*_{12}{\mathcal U}_{\xi}\otimes p^*_1E_1 \otimes
p^*_{23}{\mathcal U}^*_\xi (-{\Delta}))\, \longrightarrow\, 
{\mathcal R}^1 \, \stackrel{\alpha}{\longrightarrow}\,
E_1\otimes TX
$$
\begin{equation}\label{17}
\longrightarrow \,  {\mathcal R}^2 p_{3*}(
p^*_{12}{\mathcal U}_{\xi}\otimes p^*_1E_1 \otimes
p^*_{23}{\mathcal U}^*_\xi (-{\Delta}))
\longrightarrow\ldots
\end{equation}
For any $x\in X$ we have the cohomology exact sequence
$$
\longrightarrow\, H^i(X\times {\mathcal M}_{\xi}, \,{\mathcal U}_{\xi}
\otimes p^*_XE_1 \otimes
p^*_{{\mathcal M}_{\xi}}{\mathcal U}^*_x(-x)) \, \longrightarrow\,
 H^i(X\times {\mathcal M}_{\xi}, \, {\mathcal U}_{\xi}
\otimes p^*_XE_1 \otimes
p^*_{{\mathcal M}_{\xi}}{\mathcal U}^*_x)
$$
\begin{equation}\label{18}
\longrightarrow\, H^i({\mathcal M}_{\xi} ,\,
{\mathcal U}_x\otimes {\mathcal U}^*_x)\otimes (E_1)_x
\, \longrightarrow\,
H^{i+1}(X\times {\mathcal M}_{\xi}, \,{\mathcal U}_{\xi}
\otimes p^*_XE_1 \otimes p^*_{{\mathcal M}_{\xi}}{\mathcal U}^*_x(-x))\, ,
\end{equation}
where $(E_1)_x$ is the fibre of $E_1$ at $x$.

By \eqref{4}, $p_{X*} ({\mathcal U}_{\xi}\otimes
p_{{\mathcal M}_{\xi}}^*{\mathcal U}^*_x) =0$. So the Leray spectral
sequence for $p_X$ gives 
$$
H^1(X\times{\mathcal M}_\xi ,\,
{\mathcal U}_{\xi}\otimes p^*_XE_1\otimes 
p_{{\mathcal M}_{\xi}}^*{\mathcal U}^*_x)\,\cong\, H^0(X,\, {\mathcal
R}^1 p_{X*} ({\mathcal U}_{\xi}\otimes p^*_XE_1\otimes 
p_{{\mathcal M}_{\xi}}^*{\mathcal U}^*_x))
$$
and
$$
H^1(X\times{\mathcal M}_\xi ,\,
{\mathcal U}_{\xi}\otimes p^*_XE_1\otimes 
p_{{\mathcal M}_{\xi}}^*{\mathcal U}^*_x(-x))\,\cong\, H^0(X,\,
{\mathcal R}^1 p_{1*} ({\mathcal U}_{\xi}\otimes p^*_XE_1\otimes 
p_{{\mathcal M}_{\xi}}^*{\mathcal U}^*_x)(-x)).
$$
Since ${\mathcal U}_x$ is simple \cite[Theorem 2]{nr},
\eqref{18} gives the exact sequence 
$$
0\,\longrightarrow\, (E_1)_x\, \longrightarrow \, H^0(X,{\mathcal
R}^1 p_{X*}({\mathcal U}_{\xi}\otimes p^*_XE_1\otimes 
p_{{\mathcal M}_{\xi}}^*{\mathcal U}^*_x)(-x))\hspace{3.cm}
$$
$$
\hspace{3.cm} \longrightarrow\,
H^0(X,{\mathcal R}^1 p_{X*}
({\mathcal U}_{\xi}\otimes p^*_XE_1\otimes 
p_{{\mathcal M}_{\xi}}^*{\mathcal U}^*_x)) \, \longrightarrow\ldots
$$
This implies that ${\mathcal R}^1 p_{X*}
({\mathcal U}_{\xi}\otimes p^*_XE_1\otimes 
p_{{\mathcal M}_{\xi}}^*{\mathcal U}^*_x)$ has torsion at $x$. Now from
\eqref{16} we conclude that ${\mathcal R}^1 p_{X*}
({\mathcal U}_{\xi}\otimes p^*_XE_1\otimes 
p_{{\mathcal M}_{\xi}}^*{\mathcal U}^*_x)$ 
is a torsion sheaf, and hence
$$
H^1(X,\, {\mathcal R}^1 p_{X*}
({\mathcal U}_{\xi}\otimes p^*_XE_1\otimes 
p_{{\mathcal M}_{\xi}}^*{\mathcal U}^*_x)(-x))\, =\, H^1(X,\,
{\mathcal R}^1 p_{X*} ({\mathcal U}_{\xi}\otimes p^*_XE_1\otimes 
p_{{\mathcal M}_{\xi}}^*{\mathcal U}^*_x))\,=\,0\, .
$$

The
Leray spectral sequence for $p_X$ now yields
\begin{equation}\label{19}
H^2(X\times{\mathcal M}_\xi ,\,
{\mathcal U}_{\xi}\otimes p^*_XE_1\otimes 
p_{{\mathcal M}_{\xi}}^*{\mathcal U}^*_x)\, \cong \, H^0(X,\, {\mathcal
R}^2 p_{X*} ({\mathcal U}_{\xi}\otimes p^*_XE_1\otimes 
p_{{\mathcal M}_{\xi}}^*{\mathcal U}^*_x))
\end{equation}
and
\begin{equation}\label{20}
H^2(X\times{\mathcal M}_\xi ,\,
{\mathcal U}_{\xi}\otimes p^*_XE_1\otimes 
p_{{\mathcal M}_{\xi}}^*{\mathcal U}^*_x(-x))\cong H^0(X,{\mathcal R}^2
p_{X*} ({\mathcal U}_{\xi}\otimes p^*_XE_1\otimes 
p_{{\mathcal M}_{\xi}}^*{\mathcal U}^*_x)(-x))\, .
\end{equation}
Now from \eqref{16} it follows that 
${\mathcal R}^2 p_{X*}
({\mathcal U}_{\xi}\otimes p^*_XE_1\otimes 
p_{{\mathcal M}_{\xi}}^*{\mathcal U}^*_x)$ is a torsion sheaf, and from
Proposition \ref{3.5} and \eqref{19} that its space of sections is $0$.
So ${\mathcal R}^2 p_{X*}
({\mathcal U}_{\xi}\otimes p^*_XE_1\otimes 
p_{{\mathcal M}_{\xi}}^*{\mathcal U}^*_x)=0$ and by \eqref{20} we have
\begin{equation}\label{eqh2}
H^2(X\times{\mathcal M}_\xi ,\,
{\mathcal U}_{\xi}\otimes p^*_XE_1\otimes 
p_{{\mathcal M}_{\xi}}^*{\mathcal U}^*_x(-x))\,=\, 0\, .
\end{equation}

\begin{propn}\label{4.2}
$ {\mathcal R}^2 p_{3*}( p^*_{12}
{\mathcal U}_{\xi}\otimes p^*_1E_1 \otimes
p^*_{23}{\mathcal U}^*_\xi (-{\Delta}))\,= \,0.$
\end{propn}

\begin{proof} This follows at once from \eqref{eqh2}.\end{proof}

\begin{propn}\label{2.5}
${\mathcal R}^1 \cong E_1\otimes TX.$
\end{propn}

\begin{proof} By Corollary \ref{r1},
${\mathcal R}^1$ is a vector bundle of rank $n_0$.
Moreover Proposition \ref{4.2} implies that  the map
$\alpha \, :\, {\mathcal R}^1 \, \longrightarrow\, 
E_1\otimes TX$ in the exact sequence \eqref{17} 
is surjective. Therefore $\alpha$ must be an isomorphism,
and the proof is complete. 
\end{proof}

\section{An Inversion Formula}
We are now ready to prove our inversion formula.
\begin{thm}\label{2.7}Suppose that \eqref{deg} holds and that $E$ 
is a semistable bundle of rank $n_0$ and degree $d_0$. Then 
$$E\cong{\mathcal R}^1p_{X*}(p_{{\mathcal M}_{\xi}}^*{\mathcal W}_{\xi}(E)
\otimes{\mathcal U}_{\xi}^*\otimes p_X^*K_X).$$
\end{thm}

\begin{proof}We consider the Leray exact sequence for the composite
$$p_X\circ p_{23}=p_3:X\times{\mathcal M}_{\xi}\times X\longrightarrow X.$$
Note that ${\mathcal R}^ip_{23*}(p_{12}^*{\mathcal U}_{\xi}\otimes p_1^*E
\otimes p_{23}^*{\mathcal U}_{\xi}^*)=0$ for $i\ge1$ by \eqref{deg}, and
\begin{eqnarray*}p_{23*}(p_{12}^*{\mathcal U}_{\xi}\otimes p_1^*E
\otimes p_{23}^*{\mathcal U}_{\xi}^*)&\cong&
p_{{\mathcal M}_{\xi}}^*(p_{{\mathcal M}_{\xi}*}({\mathcal U}_{\xi}
\otimes p_X^*E))\otimes{\mathcal U}_{\xi}^*\\
&\cong&p_{{\mathcal M}_{\xi}}^*{\mathcal W}_{\xi}(E)\otimes{\mathcal U}_
{\xi}^*.
\end{eqnarray*}

Now the Leray spectral sequence gives
\begin{eqnarray*}
{\mathcal R}^1p_{X*}(p_{{\mathcal M}_{\xi}}^*
{\mathcal W}_{\xi}(E)
\otimes{\mathcal U}_{\xi}^*)&\cong&{\mathcal R}^1p_{3*}
(p_{12}^*{\mathcal U}_{\xi}\otimes p_1^*E
\otimes p_{23}^*{\mathcal U}_{\xi}^*)\\
&\cong&E\otimes TX
\end{eqnarray*}
by Proposition \ref{2.5}. Tensoring with $K_X$ now gives the result.
\end{proof}

\section{Infinitesimal Deformations}
In this section we turn to the computation of the infinitesimal deformations
of the generalised Picard bundle.
\begin{thm}\label{2.6}Suppose that \eqref{deg} holds and that $E_1$, $E_2$ 
are semistable bundles of rank $n_0$ and degree $d_0$. Then
$$H^0({\mathcal M}_{\xi},{\mathcal W}_{\xi}(E_1)\otimes 
{\mathcal W}_{\xi}(E_2)^*)\cong H^0(X,E_1\otimes E^*_2).$$ 
\end{thm}

\begin{proof}
By Propositions \ref{leray} and \ref{2.3}, we have 
$$ H^0({\mathcal M}_{\xi},\, {\mathcal W}_{\xi}(E_1)\otimes 
{\mathcal W}_{\xi}(E_2)^*) \, \cong\, H^0(X,\, {\mathcal R}^1\otimes
 E^*_2\otimes K_X)\, .
$$
{}From Proposition \ref{2.5} it follows immediately that
$$
H^0(X,\, {\mathcal R}^1\otimes E^*_2\otimes K_X)
\, \cong \, H^0(X,\, E_1\otimes E^*_2)
$$
and hence the proof is complete.
\end{proof}

\begin{coro}\label{2.8}
If $E$ is semistable and simple of rank $n_0$ and degree $d_0$, then
$$H^0({\mathcal M}_{\xi},End({\mathcal W}_{\xi}(E))
\cong {\mathbb C}.$$
In other words, the vector bundle ${\mathcal W}_{\xi}(E)$ is simple.
\end{coro}

\begin{proof}
Take $E_1=E_2=E$ in the theorem. 
\end{proof}

The following theorem now gives the infinitesimal deformations
of ${\mathcal W}_{\xi}(E)$.

\begin{thm}\label{2.9} Suppose that \eqref{deg} holds.
For any semistable bundle $E$ of rank $n_0$ and degree $d_0$,
the space of infinitesimal deformations of the vector bundle
${\mathcal W}_{\xi}(E)$, namely $H^1({\mathcal M}_\xi ,\, End
( {\mathcal W}_{\xi}(E)))$, is canonically isomorphic to
$H^1(X,\, End E)$. In particular, if $E$ is also simple, 
$$
\dim H^1({\mathcal M}_\xi ,\, End  ( {\mathcal W}_{\xi}(E)))
\, =\, n_0^2(g-1)+1.
$$
\end{thm}

\begin{proof}
Let $E_1=E_2=E$. From Propositions \ref{leray} and \ref{2.4} we
obtain an isomorphism
\begin{equation}\label{5}
H^1(X,\, {\mathcal R}^1\otimes E^* \otimes 
 K_X)\, \cong H^1({\mathcal M}_\xi ,\, End ( {\mathcal W}_{\xi}(E)))
\end{equation}

{}From Proposition \ref{2.5} we have
$$\begin{array}{lll}
H^1(X,\, {\mathcal R}^1 \otimes 
E^*\otimes K_X)\,& \cong& H^1(X,E\otimes TX \otimes E^*\otimes K_X)\\  
&\cong &H^1(X, \, End E)\, .
\end{array}
$$
Hence
$$
H^1({\mathcal M}_\xi ,\, End ( {\mathcal W}_{\xi}(E)))
\,\cong \, H^1(X, \, End E)\, 
$$
as required. The formula for the dimension follows from Riemann-Roch.
\end{proof}

\begin{rem}\label{2.10}
\begin{em} From the proof of Theorem
\ref{2.9} we see that 
$$H^1({\mathcal M}_\xi ,\, {\mathcal W}_{\xi}(E_1)
\otimes {\mathcal W}_{\xi}(E_2)^*)
\,\cong \, H^1(X, \, E_1\otimes E^*_2)\, $$
for any semistable $E_1$, $E_2$ of rank $n_0$ and degree $d_0$. 
In particular, if $E_1$, $E_2$ are stable and 
not isomorphic, we have $H^0({\mathcal M}_\xi ,\, {\mathcal W}_{\xi}(E_1)
\otimes {\mathcal W}_{\xi}(E_2)^*)=0$ by Theorem \ref{2.6} and hence
$${\dim} H^1({\mathcal M}_\xi ,\, {\mathcal W}_{\xi}(E_1)
\otimes {\mathcal W}_{\xi}(E_2)^*)\, =\, n_0^2(g-1)\, .$$
\end{em}
\end{rem}

\medskip

\section{Family of Deformations}

In this section we investigate local and global
deformations of the generalised Picard bundles constructed above.

First suppose that $(n_0,d_0)=1$ and let ${\mathcal U}(n_0,d_0)$ be a
universal bundle over $X\times{\mathcal M}(n_0,d_0)$. Now consider
$X\times{\mathcal M}(n_0,d_0)\times{\mathcal M}_{\xi}$ and define 
$$\widehat{{\mathcal U}_{\xi}}:=
p^*_{12}{\mathcal U}(n_0,d_0)\otimes p^*_{13}{\mathcal U}_{\xi},$$
and
$$\widehat{{\mathcal W}_{\xi}}:= p_{23*}(\widehat{{\mathcal U}_{\xi}})$$
By \eqref{deg} we have ${\mathcal R}^i p_{23*}
(\widehat{{\mathcal U}_{\xi}}) =0$ for
$i\not= 0$ and so $\widehat{{\mathcal W}_{\xi}}$ is locally free. 
Moreover
$$\widehat{{\mathcal W}_{\xi}}\big\vert _{\{E\}\times {\mathcal M}_{\xi}}
\cong  {\mathcal W}_{\xi}(E), $$
so $\widehat{{\mathcal W}_{\xi}}$ is a family of deformations of
${\mathcal W}_{\xi}(E)$.

\begin{thm}\label{5.1}The family $\widehat{{\mathcal W}_{\xi}}$ is
injectively parametrised and is locally complete at every point 
$E_0\in{\mathcal M}(n_0,d_0)$.
\end{thm}

\begin{proof}
The injectivity follows from Theorem \ref{2.7}. It remains to prove
that the infinitesimal deformation map is an isomorphism at every $E_0$.
By Corollary \ref{2.8}, ${\mathcal W}_{\xi}(E_0)$ is simple, so possesses 
a local analytic moduli space $S$. It follows that there exists a
neighbourhood $U$ of $E_0$ in ${\mathcal M}(n_0,d_0)$ with respect to the
analytic topology and a holomorphic map
$$\phi:U\longrightarrow S$$
such that $\phi(E)\cong{\mathcal W}_{\xi}(E)$ for all $E\in U$.
By injectivity the image of $\phi$ has dimension
$$\dim{\mathcal M}(n_0,d_0)=n_0^2(g-1)+1$$
at every point. On the other hand, by Theorem \ref{2.9}, we know that the
Zariski tangent space to $S$ at ${\mathcal W}_{\xi}(E)$ also has dimension
$n_0^2(g-1)+1$. It follows that $S$ is smooth at ${\mathcal W}_{\xi}(E)$.
Hence, by Zariski's Main Theorem, $\phi$ maps $U$ isomorphically onto an
open subset of $S$, and in particular the differential $d\phi$ (which
coincides with the infinitesimal deformation map) is an isomorphism at $E_0$.
\end{proof}

\begin{rem}\label{5.2}
\begin{em}When $(n_0,d_0)\ne1$, we no longer have a universal
bundle ${\mathcal U}(n_0,d_0)$. However, for any $E\in
{\mathcal M}(n_0,d_0)$,
there exists an \'etale neighbourhood
of $E_0$ over which a universal bundle does exist. The argument of 
Theorem \ref{5.1} now goes through to give a family of Picard bundles
which is locally complete at $E_0$. This family is not injectively 
parametrised, but it is still true that
$${\mathcal W}_{\xi}(E_1)\cong{\mathcal W}_{\xi}(E_2)\Leftrightarrow
E_1\cong E_2.$$
\end{em}\end{rem}

Theorem \ref{5.1} says that ${\mathcal M}(n_0,d_0)$ 
is in some sense a
moduli space for the bundles ${\mathcal W}_{\xi}(E)$. Since 
${\mathcal M}(n_0,d_0)$ is irreducible, this implies that all
${\mathcal W}_{\xi}(E)$ have the same Hilbert polynomial $P_{n_0,d_0}$
with respect to the unique polarisation $\theta_{\xi}$ of 
${\mathcal M}_{\xi}$. We do not know in general that the 
${\mathcal W}_{\xi}(E)$ 
possess a good global moduli space. However, if all ${\mathcal W}_{\xi}(E)$
are stable with respect to $\theta_{\xi}$, then they belong to the
moduli space ${\mathcal M}_{{\mathcal M}_{\xi}}(P_{n_0,d_0})$ and indeed to
one particular connected component ${\mathcal M}^0$ of this moduli space.
The map $E\mapsto{\mathcal W}_{\xi}(E)$ then defines a morphism
$$\phi:{\mathcal M}(n_0,d_0)\longrightarrow{\mathcal M}^0.$$

\begin{thm}\label{5.3}If $(n_0,d_0)=1$ and ${\mathcal W}_{\xi}(E)$ is 
stable with
respect to $\theta_{\xi}$ for every $E\in{\mathcal M}(n_0,d_0)$, then $\phi$ 
is an isomorphism.
\end{thm}

\begin{proof}By Theorem \ref{5.1}, $\phi$ is an isomorphism onto an
open subset of ${\mathcal M}^0$. Since ${\mathcal M}(n_0,d_0)$ is
complete, this implies that $\phi$ is an isomorphism.
\end{proof}

Theorem \ref{5.3} applies in particular if $n_0=1$. In this case we can
suppose that $d_0=0$, so that ${\mathcal M}(n_0,d_0)=J$, the Jacobian of $X$.
We know by \cite{bbgn} that ${\mathcal W}_{\xi}({\mathcal O})$ is stable with
respect to $\theta_{\xi}$, and the same proof shows that
${\mathcal W}_{\xi}(L)$ is stable for any $L\in J$.

In this case, we can go a little further. Since we shall want to allow
$X$ and $\xi$ to vary, we denote the
space ${\mathcal M}^0$ by ${\mathcal M}^0_{X,\xi}$. Let $\Theta$ denote the 
principal polarisation on $J$ defined by a theta divisor and let $\zeta$ 
denote the  polarisation on ${\mathcal M}^0_{X,\xi}$ defined
by the determinant line bundle \cite[Section 4]{ib}.

\begin{thm}\label{5.4}
With respect to the above polarisations, the morphism
$$
\phi\, : \,J\, \longrightarrow\,
{\mathcal M}^0_{X,\xi}
$$ 
is an isomorphism of polarised varieties.
\end{thm}

\begin{proof}
We wish to show that the isomorphism $\phi$ takes
$\zeta$ to a nonzero constant scalar multiple (independent of the
curve $X$) of $\Theta$.

Take any family of pairs $(X,\xi)$, where $X$ is a connected
non-singular projective curve of genus
$g$ and $\xi$ is a line bundle on $X$ of degree $d>n(2g-2)$,
parametrized  by a connected space
$T$. Consider the corresponding family of moduli spaces ${\mathcal
M}^0_{X,\xi}$ (respectively, Jacobians $J$) over $T$,
where
$(X,\xi)$ runs over the family. Using the map $\phi$ an isomorphism
between these two families is obtained.
The polarisation $\zeta$ (respectively, $\Theta$) defines a constant
section of the second direct image over $T$
of the constant sheaf $\mathbb Z$ over the family.
It is known that for the general curve $X$ of genus $g$, the
Neron-Severi group of $J$ is $\mathbb Z$. Therefore,
for such a curve, $\phi$ takes $\zeta$ to a nonzero constant scalar
multiple of $\Theta$. Since $T$ is connected,
if $T$ contains a curve with
$NS(J)\, =\, {\mathbb Z}$, then $\phi$
takes $\zeta$ to the same nonzero constant scalar multiple
of $\Theta$ for every curve in the family. Since the moduli
space of smooth curves of genus $g$ is connected, the proof
is complete.
\end{proof}

Finally we have our Torelli theorem. 

\begin{coro}\label{5.5}
Let $X$ and $X'$ be two non-singular algebraic curves of 
genus $g\geq 3$ and let $\xi$ (respectively $\xi'$) be a line bundle
of degree $d>n(2g-2)$ on $X$ (respectively $X'$). If 
${\mathcal M}^0_{X,\xi}\cong 
 {\mathcal M}^0_{X',\xi'}$ as polarised varieties
then 
$X\cong X'$.
\end{coro}

\begin{proof}
This follows at once from Theorem \ref{5.4} and the
classical Torelli theorem.
\end{proof}


\end{document}